\documentclass[11pt]{article}
\usepackage{latexsym}
\usepackage[all]{xy}
\usepackage{amsmath}
\usepackage{amssymb}
\usepackage{amsfonts}

\newtheorem{thm}{Theorem}[section]
\newtheorem{prop}[thm]{Proposition}
\newtheorem{lem}[thm]{Lemma}

\newtheorem{df}[thm]{Definition}
\newtheorem{cor}[thm]{Corollary}

\newtheorem{rmk}[thm]{Remark}

\begin{document}

\title{\textbf{Derived Hall algebras}}
\bigskip
\bigskip

\author{\bigskip\\
Bertrand To\"en \\
\small{Laboratoire Emile Picard}\\
\small{UMR CNRS 5580} \\
\small{Universit\'{e} Paul Sabatier, Toulouse}\\
\small{France}}

\date{June 2005}

\maketitle

\begin{abstract}
The purpose of this work is to define a derived Hall algebra
$\mathcal{DH}(T)$, associated to any dg-category $T$ (under
some finiteness conditions), generalizing the Hall algebra
of an abelian category. Our main theorem states that
$\mathcal{DH}(T)$ is associative and unital. When the
associated triangulated category $[T]$ is endowed
with a t-structure with heart $\mathcal{A}$, it is shown that 
$\mathcal{DH}(T)$ contains the usual Hall algebra
$\mathcal{H}(\mathcal{A})$.
We will also prove
an explicit formula for the derived Hall numbers purely in terms of
invariants of the triangulated category associated to $T$. As an example, 
we describe the derived Hall algebra of an hereditary abelian category.
\end{abstract}

\medskip


\tableofcontents

\bigskip

\section{Introduction}

Hall algebras are associative algebras associated to
abelian categories (under some finiteness conditions).
They appear, as well as various variants,
in several mathematical contexts, in
particular in the theory of quantum groups, and
have been intensively studied (we refer to
the survey paper \cite{de} for an overview of the
subject as well as for bibliographic references).
Roughly speaking, the Hall algebra $\mathcal{H}(A)$ of
an abelian category $A$, is generated as a vector
space by the isomorphism classes of objects in $A$, and
the product  of two objects $X$ and $Y$ is defined
as  $\sum_{Z\in Iso(A)}g_{X,Y}^{Z}Z$, where
the structure constants $g_{X,Y}^{Z}$ count
a certain number of extensions $X\rightarrow Z \rightarrow Y$.
If $\mathfrak{g}$ is a  complex Lie algebra with
simply laced Dynkin diagram,
and $A$ the abelian category of $\mathbb{F}_{q}$-representations of
a quiver whose underlying graph is the Dynkin diagram
of $\mathfrak{g}$, Ringel proved that the Hall algebra $\mathcal{H}(A)$ is closely
related to $U_{q}(b^{+})$, the positive part of the quantum enveloping algebra $U_{q}(\mathfrak{g})$
(see \cite{ri}).
It seems to have been stressed by several authors that
an extension of the construction $A\mapsto \mathcal{H}(A)$
to the derived category of $A$-modules could be useful to
recover the full algebra $U_{q}(g)$. To this end, certain
Hall numbers $g_{X,Y}^{Z}$, that are supposed to count
triangles $\xymatrix{X \ar[r] & Z \ar[r] & Y \ar[r] & X[1]}$
in a triangulated category, have been introduced
(see for example \cite{de} \S 2.2). However, these Hall numbers
do not provide an associative multiplication. Quoting
the introduction of \cite{ka}: \\

\emph{Unfortunately, a direct mimicking
of the Hall algebra construction, but with exact triangles
replacing exact sequences, fails to give an associative
mutiplication, even though the octohedral
axiom looks like the right tool to establish associativity.} \\

The purpose of this work is to define analogs of
Hall algebras and Hall numbers, called \emph{derived Hall algebra} and
\emph{derived Hall numbers}, associated
with dg-categories $T$ (instead of abelian or triangulated categories)
and providing an associative and unital multiplication on the
vector space generated by equivalence classes of
perfect $T^{op}$-dg-modules (see \cite{to} for definitions and notations).
Our main theorem is the following.

\begin{thm}\label{ti}
Let $T$ be a dg-category over some finite field
$\mathbb{F}_{q}$, such that for any two objects
$x$ and $y$ in $T$ the graded vector space
$H^{*}(T(x,y))=Ext^{*}(x,y)$ is of finite total dimension.
Let $\mathcal{DH}(T)$ be the $\mathbb{Q}$-vector space with basis
$\{\chi_{x}\}$, where $x$ runs through the set of equivalence classes of perfect $T^{op}$-dg-modules $x$.
Then, there
exists an associative and unital product
$$\mu : \mathcal{DH}(T) \otimes \mathcal{DH}(T) \longrightarrow
\mathcal{DH}(T),$$
with
$$\mu(\chi_{x},\chi_{y})=\sum_{z}g_{x,y}^{z}\chi_{z},$$
and where the $g_{x,y}^{z}$ are rational numbers counting
the virtual number of morphisms in $T$, $f : x'\rightarrow z$, such that
$x'$ is equivalent to $x$ and the homotopy cofiber of $f$ is equivalent
to $y$.
\end{thm}

Our construction goes as follows. To each dg-category $T$, satisfying the
condition of the theorem, we associate a certain homotopy type $X^{(0)}(T)$,
which is a classifying space of perfect $T^{op}$-dg-modules.
The connected components of $X^{(0)}(T)$ are in one to one correspondence
with the quasi-isomorphism classes of perfect $T^{op}$-dg-modules. Furthermore,
for any point $x$, corresponding to a $T^{op}$-dg-module, the homotopy groups
of $X^{(0)}(T)$ at $x$ are given by
$$\pi_{1}(X^{(0)}(T),x)=Aut(x) \qquad \pi_{i}(X^{(0)}(T),x)=Ext^{1-i}(x,x)
\; \forall i>1$$
where the automorphism and Ext groups are computed in the
homotopy category of $T^{op}$-dg-modules. In the same way, we consider
the simplicial set $X^{(1)}(T)$ which is a classifying space for
morphisms between perfect $T^{op}$-dg-modules. We then construct
a diagram of simplicial sets
$$\xymatrix{X^{(1)}(T) \ar[r]^-{c} \ar[d]_-{s\times t} & X^{(0)}(T) \\
X^{(0)}(T)\times X^{(0)}(T), & }$$
where $c$ sends a morphism of $T$ to its target, and $s\times t$ sends
a morphism $u : x\rightarrow y$ of $T$ to the pair consisting of
$x$ and the homotopy cofiber (or cone) of $u$. We note that
the vector space $\mathcal{DH}(T)$, with basis $\{\chi_{x}\}$, can be
identified with the set of $\mathbb{Q}$-valued functions
with finite support on $\pi_{0}(X^{(0)}(T))$. Our multiplication
is then defined as a certain convolution product
$$\mu:=c_{!}\circ (s\times t)^{*} : \mathcal{DH}(T)\otimes \mathcal{DH}(T)
\longrightarrow \mathcal{DH}(T).$$
We should mention that though the operation $(s\times t)^{*}$ is the naive pull-back of
functions, the morphism $c_{!}$ is a fancy push-forward, defined in terms of the
full homotopy types $X^{(0)}(T)$ and $X^{(1)}(T)$ and involving alternating
products of the orders of their
higher homotopy groups. The fact that this product is associative is then
proved using a base change formula, and the fact that certain
commutative squares involving $X^{(2)}(T)$, the classifying space of
pairs of composable morphisms in $T$, are homotopy cartesian.

In the second part of the paper we will establish an explicit formula for the
derived Hall numbers, purely in terms of
the triangulated category of perfect dg-modules (see Prop. \ref{p3}).
This will imply an invariance statement for our
derived Hall algebras under derived equivalences. This also implies the
\emph{heart property}, stating that for
a given $t$-structure on the category of perfect $T^{op}$-dg-modules, the
sub-space of $\mathcal{DH}(T)$ spanned by the objects belonging to the 
heart $\mathcal{C}$
is a sub-algebra isomorphic to the usual Hall algebra of $\mathcal{C}$ (see Prop. \ref{p4}). 
When $A$ is a finite dimensional algebra, considered
as a dg-category $BA$, this allows us to identify the usual Hall algebra
$\mathcal{H}(A)$ as the sub-algebra in $\mathcal{DH}(BA^{op})$ spanned by the
non-dg $A$-modules. We will finish this work with a description 
by generators and relations
of the derived Hall algebra of an abelian category of cohomological
dimension $\leq 1$, in terms of its underived Hall algebra (see Prop. \ref{p5}). \\

\textbf{Related and future works:} As we mentioned, there have been a huge
amount of works on Hall algebras, in which several variants of them appear
(e.g. Ringel-Hall algebras). There is no
doubt that one can also define derived versions of these variant notions,
and studying these seems to me an important
work to do, as this might throw new lights on quantum groups. However, the precise relations
between derived Hall algebras and quantum groups are still unclear and require
to be investivated in the future. In fact, it seems to be expected that the
full quantum enveloping algebra is related not to the derived category of representations
of the associated quiver but rather to its 2-periodic version. A 2-periodic 
dg-category will never satisfy the finiteness condition needed in this work in order to construct
the derived Hall algebra, and therefore it seems that though we think the present work
is a step in the right direction it is probably not enough to answer the question of
reconstructing the full quantum group as the Hall algebra of something. 

G. Lusztig has given a very nice geometric construction of Hall algebras, using
certain perverse sheaves on the moduli stacks of modules over associative
algebras (see for example \cite{lu}). This approach gives a much more flexible theory, as for example it
allows to replace the finite field by essentially any base, and has found
beautiful applications to quantum groups (such as the construction of the
so-called canonical bases). I am convinced that the
derived Hall algebra $\mathcal{DH}(T)$ of a dg-category $T$ also
has a geometric interpretation in terms of complexes of l-adic sheaves
on the moduli stack of perfect $T^{op}$-dg-modules
(as constructed in \cite{to-va}). 
This geometric interpretation of the derived Hall algebras will be
investigated in a future work. \\

\textbf{Acknowledgments:} I am very grateful to B. Keller for pointing out to me the
problem of defining Hall algebras for dg-categories, and for having noticed its relation
to the work \cite{to-va}. 

I also would like to thanks both referees for their comments. \\

\bigskip

\textbf{Conventions and notations:} We use model category theory from
\cite{ho}. For a diagram $\xymatrix{x\ar[r] & z & \ar[l] y}$ in a model category $M$ we denote by
$x\times_{z}^{h}y$ the homotopy fiber product, and
by $x\coprod_{z}^{h}y$ the homotopy push-out.
For dg-categories, we refer to \cite{tab,to}, and we will mainly
keep the same notation and terminology as in \cite{to}. Finally, for a finite
set $X$ we denote by $|X|$ the cardinality of $X$.

\section{Functions on locally finite homotopy types}

In this section, we will discuss some finiteness conditions on
homotopy types, functions with
finite supports on them and their transformations (push-forwards,
pull-backs and base change).

\begin{df}\label{d1}
A homotopy type $X\in Ho(Top)$
is \emph{locally finite}, if it satisfies the following two
conditions.
\begin{enumerate}
\item For all base points $x\in X$, the
group $\pi_{i}(X,x)$ is finite.
\item For all base point $x \in X$,
there is $n$ (depending on $x$) such that
$\pi_{i}(X,x)=0$ for all $i>n$.
\end{enumerate}
\end{df}

The full sub-category of $Ho(Top)$ consisting of all
locally finite objects will be denoted by
$Ho(Top)^{lf}$. \\

\begin{df}\label{d2}
Let  $X\in Ho(Top)$. The
$\mathbb{Q}$-vector space of \emph{rational
functions with finite support on $X$} is
the $\mathbb{Q}$-vector space of functions
on the set $\pi_{0}(X)$ with values in $\mathbb{Q}$ and with finite support.
It is denoted by $\mathbb{Q}_{c}(X)$.
\end{df}

Let $f : X \longrightarrow Y$ be a morphism in $Ho(Top)^{lf}$. We define
a morphism
$$f_{!} : \mathbb{Q}_{c}(X) \longrightarrow \mathbb{Q}_{c}(Y)$$
in the following way. For $y\in \pi_{0}(Y)$, one sets
$F_{y} \in Ho(Top)$ the homotopy fiber of $f$ at the point $y$, and
$i : F_{y} \longrightarrow X$ the natural morphism.
We note that the long exact sequence in homotopy implies
that for all $z\in \pi_{0}(F_{y})$
the group $\pi_{i}(F_{y},z)$ is finite and moreover zero for
all $i$ big enough. Furthermore, the fibers of the map
$i : \pi_{0}(F_{y}) \longrightarrow \pi_{0}(X)$ are finite.

For a function $\alpha \in \mathbb{Q}_{c}(X)$,
and $y\in \pi_{0}(Y)$, one sets
$$
f_{!}(\alpha)(y)=\sum_{z\in \pi_{0}(F_{y})}
\alpha(i(z)).\prod_{i>0}|\pi_{i}(F_{y},z)|^{(-1)^{i}}.
$$
Note that as $\alpha$ has  finite support, the sum
in the definition of $f_{!}$ only has a finite number of non-zero
terms, and thus is well defined.

\begin{lem}\label{l1}
For any morphism $f : X \longrightarrow Y$ in $Ho(Top)^{lf}$,
any $\alpha\in \mathbb{Q}_{c}(X)$, and any $y\in \pi_{0}(Y)$, one has
$$f_{!}(\alpha)(y)=\sum_{x\in \pi_{0}(X) \diagup f(x)=y}
\alpha(x).\prod_{i>0}\left( |\pi_{i}(X,x)|^{(-1)^{i}}.|\pi_{i}(Y,y)|^{(-1)^{i+1}}
\right)
$$
\end{lem}

\textit{Proof:} Let $y\in \pi_{0}(Y)$, $i : F_{y} \rightarrow X$ the
homotopy fiber of $f$ at $y$, and $z\in \pi_{0}(F_{y})$.
The long exact sequence in homotopy
$$\xymatrix{
\pi_{i+1}(Y,y) \ar[r] & \pi_{i}(F_{y},z) \ar[r] &
\pi_{i}(X,x) \ar[r] & \pi_{i}(Y,y) \ar[r] &
\pi_{i-1}(F_{y},z) \ar[r] & \dots }$$
$$\xymatrix{ \dots \ar[r] & \pi_{1}(Y,y) \ar[r] &
\pi_{0}(F_{y},z) \ar[r] & \pi_{0}(X,x) \ar[r] & \pi_{0}(Y,y),}$$
implies that one has for any $x\in \pi_{0}(X)$ with $f(x)=y$
$$\prod_{i>0}\left( |\pi_{i}(X,x)|^{(-1)^{i}}.|\pi_{i}(Y,y)|^{(-1)^{i+1}}
\right)
=\sum_{z\in \pi_{0}(F_{y})\diagup i(z)=x}
\prod_{i>0}|\pi_{i}(F_{y},z)|^{(-1)^{i}}.$$
\hfill $\Box$ \\

\begin{cor}\label{cl1}
For any morphisms $f : X \longrightarrow Y$ and
$g : Y \longrightarrow Z$ in $Ho(Top)^{lf}$, one has
$(g\circ f)_{!}=g_{!}\circ f_{!}$.
\end{cor}

\textit{Proof:} Clearly follows from the lemma \ref{l1}. \hfill $\Box$ \\

\begin{df}\label{d3}
A morphism $f : X \longrightarrow Y$ in $Ho(Top)^{lf}$
is called \emph{proper},
if for any $y\in \pi_{0}(Y)$ there is only a finite
number of elements $x \in \pi_{0}(X)$ with $f(x)=y$.
\end{df}

Note that a morphism $f : X \longrightarrow Y$ in $Ho(Top)^{lf}$
is  proper if and only if for any $y\in \pi_{0}(Y)$
the homotopy fiber $F_{y}$ of $f$ at $y$ is such that
$\pi_{0}(F_{y})$ is a finite set. \\

Let $f : X \longrightarrow Y$ be a proper morphism
between locally finite homotopy types. We define
the pull-back
$$f^{*} : \mathbb{Q}_{c}(Y) \longrightarrow \mathbb{Q}_{c}(X)$$
by the usual formula
$$f^{*}(\alpha)(x)=\alpha(f(x)),$$
for any $\alpha\in \mathbb{Q}_{c}(Y)$ and any
$x\in \pi_{0}(X)$. Note that as the morphism $f$ is proper
$f^{*}(\alpha)$ has finite support, and thus that
$f^{*}$ is well defined.

\begin{lem}\label{l3}
Let
$$\xymatrix{
X' \ar[r]^-{u'} \ar[d]_-{f'} & X \ar[d]^-{f} \\
Y' \ar[r]_-{u} & Y}$$
be a homotopy cartesian square of locally finite homotopy types,
with $u$ being proper. Then $u'$ is proper and furthermore one has
$$u^{*}\circ f_{!}=(f')_{!}\circ (u')^{*} :
\mathbb{Q}_{c}(X) \longrightarrow \mathbb{Q}_{c}(Y').$$
\end{lem}

\textit{Proof:} The fact that $u'$ is proper if $u$ is so follows
from the fact that the homotopy fiber of $u'$ at a point
$x\in \pi_{0}(X)$ is equivalent to the homotopy fiber
of $u$ at the point $f(x)$.

Let $y'\in \pi_{0}(Y')$, $y=u(y')$. Let us denote by
$F$ the homotopy fiber of $f'$ at $y'$. Note that
$F$ is also the homotopy fiber of $f$ at $y$. We denote
by $i' : F \longrightarrow X'$ the natural morphism, and
we let $i=u'\circ i' : F \longrightarrow X$. For
$\alpha\in \mathbb{Q}_{c}(X)$, one has
$$(u^{*}f_{!})(\alpha)(y')=f_{!}(\alpha)(y)=
\sum_{z\in \pi_{0}(F)}\alpha(i(z)).\prod_{i>0}|\pi_{i}(F,z)|^{(-1)^{i}}$$
$$=\sum_{z\in\pi_{0}(F)}
(u')^{*}(\alpha)(i'(z)).\prod_{i>0}|\pi_{i}(F,z)|^{(-1)^{i}}
=(f')_{!}(u'^{*}(\alpha))(y').$$
\hfill $\Box$ \\

\section{Hall numbers and Hall algebras for dg-categories}

Let $k=\mathbb{F}_{q}$ be a finite field with $q$ elements.
We refer to \cite{tab,to} for the notion of dg-categories over $k$,
as well as their homotopy theory.

\begin{df}\label{d4}
A (small) dg-category $T$ over $k$ is \emph{locally finite}
if for any two objects $x$ and $y$ in $T$, the complex
$T(x,y)$ is cohomologically bounded and with finite
dimensional cohomology groups.
\end{df}

Let $T$ be a locally finite dg-category, and
let us consider $M(T):=T^{op}-Mod$, the model category
of dg-$T^{op}$-modules. This is a stable model category,
whose homotopy category $Ho(M(T))$ is therefore
endowed with a natural triangulated structure.
We will denote by $[-,-]$ the Hom's in this triangulated
category.
Recall that an object
of $M(T)$ is perfect if it belongs to the smallest
sub-category of $Ho(M(T))$ containing the
quasi-representable modules and which is stable
by retracts, homotopy push-outs and homotopy pull-backs.
Note that the perfect objects in
$M(T)$ are also the compact objects in
the triangulated category $Ho(M(T))$. The full sub-category
of perfect objects in $M(T)$ will be denoted by $P(T)$.
It is easy to see, as $T$ is locally finite, that
for any two objects $x$ and $y$ in $P(T)$,
the $k$-vector spaces $[x,y[i]]=:Ext^{i}(x,y)$
are finite dimensional and are non-zero only for a finite
number of indices $i\in \mathbf{Z}$.

We let $I$ be the category with two objects $0$ and $1$ and
a unique morphism $0 \rightarrow 1$, and we consider the
category $M(T)^{I}$, of morphisms in $M(T)$, endowed with
its projective model structure (for which equivalences and
fibrations are defined levelwise). Note that cofibrant
objects in $M(T)^{I}$ are precisely the morphisms
$u : x\rightarrow y$, where $x$ and $y$ are cofibrant
and $u$ is a cofibration in $M(T)$.

There exist three left Quillen functors
$$s,c,t : M(T)^{I} \longrightarrow M(T)$$
defined, for any $u : x \rightarrow y$ object in $M(T)^{I}$, by
$$s(u)=x \qquad c(u)=y \qquad t(u)=y\coprod_{x}0.$$
We now consider the diagram of left Quillen functors
$$\xymatrix{
M(T)^{I} \ar[r]^-{c} \ar[d]_-{s\times t} & M(T) \\
M(T)\times M(T).}$$
Restricting to the sub-categories of cofibrant
and perfect objects, and of equivalences between them, one gets a
diagram of categories and functors
$$\xymatrix{
w(P(T)^{I})^{cof} \ar[r]^-{c} \ar[d]_-{s\times t} & wP(T)^{cof} \\
wP(T)^{cof}\times wP(T)^{cof}.}$$
Applying the nerve construction, one gets
a diagram of homotopy types
$$\xymatrix{
N(w(P(T)^{I})^{cof}) \ar[r]^-{c} \ar[d]_-{s\times t} & N(wP(T)^{cof}) \\
N(wP(T)^{cof})\times N(wP(T)^{cof}).}$$
This diagram will be denoted by
$$\xymatrix{
X^{(1)}(T) \ar[r]^-{c} \ar[d]_-{s\times t} & X^{(0)}(T) \\
X^{(0)}(T)\times X^{(0)}(T). & }$$

\begin{lem}\label{l4}
\begin{enumerate}
\item The homotopy types
$X^{(1)}(T)$ and $X^{(0)}(T)$ are locally finite.
\item
The morphism $s\times t$ above is proper.
\end{enumerate}
\end{lem}

\textit{Proof:} (1) For $x \in \pi_{0}(X^{(0)}(T))$,
corresponding to a perfect $T^{op}$-module, we know by
\cite{dk} (see also \cite[Appendix A]{hagII}) that $\pi_{1}(X^{(0)}(T),x)$
is a subset in $Ext^{0}(x,x)=[x,x]$, and that
$\pi_{i}(X^{(0)}(T),x)=Ext^{1-i}(x,x)$ for all $i>1$.
By assumption on $T$ we know that these Ext groups are finite
and zero for $i$ big enough.  This implies that
$X^{(0)}(T)$ is locally finite. The same is true
for $X^{(1)}(T)$, as it can be identified,
up to equivalence, with $X^{(0)}(T')$, for
$T'=T\otimes_{k}I_{k}$ (where $I_{k}$ is the
dg-category freely generated by the category $I$). \\

(2) Let $(x,z)$ be a point in $X^{(0)}(T)\times X^{(0)}(T)$, and
let $F$ be the homotopy fiber of $s\times t$ at $(x,z)$. It
is not hard to check that
$\pi_{0}(F)$ is in bijection with $Ext^{1}(z,x)$, and thus is finite
by assumption on $T$. This shows that $s\times t$ is proper. \hfill $\Box$ \\

The space of functions with finite support
on $X^{(0)}(T)$ will be denoted by
$$\mathcal{DH}(T):=\mathbb{Q}_{c}(X^{(0)}(T)).$$
Note that there exists a natural isomorphism
$$\begin{array}{ccc}
\mathcal{DH}(T)\otimes \mathcal{DH}(T) & \simeq&
\mathbb{Q}_{c}(X^{(0)}(T)\times X^{(0)}(T)) \\
f\otimes g & \mapsto & \left(   (x,y) \mapsto f(x).g(y) \right)  
\end{array}.$$
Therefore, using our lemma \ref{l4} one gets a natural
morphism
$$\mu:=c_{!}\circ (s\times t)^{*} : \mathcal{DH}(T)\otimes
\mathcal{DH}(T) \longrightarrow \mathcal{DH}(T).$$

Let us now consider the set
$\pi_{0}(X^{(0)}(T))$. This set is clearly in natural bijection
with the set of isomorphism classes of perfect
objects in $Ho(M(T))$. Therefore, the $\mathbb{Q}$-vector space
$\mathcal{DH}(T)$ is naturally isomorphic to the
free vector space over the isomorphism classes of perfect
objects in $Ho(M(T))$.
For an isomorphism class of a perfect $T^{op}$-module $x$, we simply denote by
$\chi_{x} \in \mathcal{DH}(T)$ the corresponding characteristic function.
When $x$ runs throught the set of isomorphism classes of perfect
objects in $Ho(M(T))$, the elements $\chi_{x}$ form a $\mathbb{Q}$-basis of $\mathcal{DH}(T)$.

\begin{df}\label{d5}
Let $T$ be a locally finite dg-category over $k$.
\begin{enumerate}
\item
The \emph{derived Hall algebra of $T$} is
$\mathcal{DH}(T)$ with the multiplication
$\mu$ defined above.
\item Let $x$, $y$ and $z$ be three perfect $T^{op}$-modules.
The \emph{derived Hall number} of $x$ and $y$ along $z$ is defined by
$$g_{x,y}^{z}:=\mu(\chi_{x},\chi_{y})(z),$$
where $\chi_{x}$ and $\chi_{y}$ are the characteristic functions
of $x$ and $y$.
\end{enumerate}
\end{df}

\section{Associativity and unit}

The main theorem of this work is the following.

\begin{thm}\label{t1}
Let $T$ be a locally finite dg-category. Then the
derived Hall algebra $\mathcal{DH}(T)$ is associative and unital.
\end{thm}

\textit{Proof:} Let us start with the existence of a unit. Let $\mathbf{1}\in
\mathcal{DH}(T)$ be the characteristic function of the zero $T^{op}$-module.
We claim that $\mathbf{1}$ is a unit for the multiplication $\mu$.
Indeed, let $x \in \pi_{0}(X^{(0)}(T))$ and $\chi_{x}\in \mathcal{DH}(T)$ its characteristic function.
The set $\pi_{0}(X^{(1)}(T))$
is identified with the isomorphism classes of perfect objects in $Ho(M(T)^{I})$
(i.e. $\pi_{0}(X^{(1)}(T))$ is the set of equivalence classes of morphisms
between perfect $T^{op}$-modules). Let us fix a morphism $u : y \rightarrow z$
between perfect $T^{op}$-modules, considered as an element in $\pi_{0}(X^{(1)}(T))$.
The function $(s\times t)^{*}(\mathbf{1},\chi_{x})$
sends $u$ to $1$ if $y$ is isomorphic to zero in $Ho(M(T))$ and $x$ is isomorphic
to $z$ in $Ho(M(T))$, and to $0$ otherwise. In other words,
$(s\times t)^{*}(\mathbf{1},\chi_{x})$ is the characteristic function of the
sub-set of $\pi_{0}(X^{(1)}(T))$ consisting of morphisms $0 \rightarrow z$ with
$z\simeq x$ in $Ho(M(T))$. We let $X$ be the sub-simplicial set of
$X^{(1)}(T)$ which is the union of all connected component
belongings to the support of $(s\times t)^{*}(\mathbf{1},\chi_{x})$.
This simplicial set is connected, and by definition of the product $\mu$, one has (see Lemma \ref{l1})
$$\mu(\mathbf{1},\chi_{x})(x)=\prod_{i>0}\left( |\pi_{i}(X)|^{(-1)^{i}}.|\pi_{i}(X^{(0)}(T),x)|^{(-1)^{i+1}}
\right) .$$
$$\mu(\mathbf{1},\chi_{x})(y)=0 \; if \; y\neq x.$$
The left Quillen functor $c : M(T)^{I} \longrightarrow M$ is
clearly fully faithful up to homotopy when restricted to the sub-category of morphisms
$y\rightarrow z$ such that $y$ is equivalent to $0$. This implies that the induced
morphism
$$c : X \longrightarrow X^{(0)}(T)$$
induces an isomorphism on all homotopy groups $\pi_{i}$ for $i>0$, identifying $X$ with
a connected component of $X^{(0)}(T)$. Therefore, one has
$\mu(\mathbf{1},\chi_{x})(x)=1$, showing that
$\mu(\mathbf{1},\chi_{x})=\chi_{x}$. Essentially the same argument shows that
$\mu(\chi_{x},\mathbf{1})=\chi_{x}$, proving that $\mathbf{1}$ is a unit for
$\mu$. \\

We now pass to the associativity of the multiplication. For sake of simplicity we first introduce
several notations. Let
$M(T)^{(1)}$ be the model category $M(T)^{I}$. We also let
$M(T)^{(2)}$ be the model category of composable pairs
of morphisms in $M$ (i.e. objects in $M(T)^{(2)}$ are
diagrams of objects in $M(T)$, $x\rightarrow y\rightarrow z$), again endowed
with its projective model structure. The full sub-categories of perfect
and cofibrant objects in $M(T)^{(i)}$ will be denoted by $(P(T)^{(i)})^{cof}$ (for $i=1,2$).
Finally, the nerve of the sub-category of equivalences will simply be denoted by
$$X^{(i)}:=N(w(P(T)^{(i)})^{cof}) \qquad i=1,2$$
$$X^{(0)}:=N(wP(T)^{cof}).$$
We recall the existence of the three Quillen functors
$$s,t,c : M(T)^{(1)} \longrightarrow M(T)$$
and their induced morphisms on the corresponding locally finite homotopy types
$$s,t,c : X^{(1)} \longrightarrow X^{(0)}.$$
We introduce four left Quillen functors
$$f,g,h,k : M(T)^{(2)} \longrightarrow M(T)^{(1)}$$
defined, for an object $u=(x\rightarrow y\rightarrow z)$ in $M(T)^{(2)}$ by
$$f(u)=x\rightarrow y \qquad g(u)=y\rightarrow z \qquad \qquad h(u)=x \rightarrow z
\qquad k(u)=y\coprod_{x}0 \rightarrow Z\coprod_{x}0.$$
These left Quillen functors induce four morphisms on the corresponding locally finite homotopy types
$$f,g,h,k : X^{(2)} \longrightarrow X^{(1)}.$$
We know consider the following two diagrams of simplicial sets
$$\xymatrix{
X^{(2)} \ar[r]^-{g} \ar[d]_-{f\times (t\circ k)} & X^{(1)} \ar[r]^-{c} \ar[d]^-{s\times t} & X^{(0)} \\
X^{(1)}\times X^{(0)} \ar[r]^-{c\times id} \ar[d]_-{(s\times t)\times id} & X^{(0)}\times X^{(0)} & \\
(X^{(0)} \times X^{(0)}) \times X^{(0)} & & \\}$$

$$\xymatrix{
X^{(2)} \ar[r]^-{h} \ar[d]_-{(s\circ f)\times k} & X^{(1)} \ar[r]^-{c} \ar[d]^-{s\times t} & X^{(0)} \\
X^{(0)}\times X^{(1)} \ar[r]^-{id\times c} \ar[d]_-{id\times (s\times t)} & X^{(0)}\times X^{(0)} & \\
X^{(0)} \times (X^{(0)} \times X^{(0)}). & & \\}$$

The left and top sides of these two diagrams are clearly equal as diagrams in $Ho(Top)^{lf}$ to
$$\xymatrix{
X^{(2)} \ar[r]^-{c\circ h} \ar[d]_-{(s\circ f)\times (s\circ k)\times (t\circ h)} & X^{(0)} \\
X^{(0)}. & \\}$$

Therefore, lemma \ref{l3} would imply that
$$c_{!}\circ (s\times t)^{*}\circ (c\times id)_{!} \circ ((s\times t)\times id)^{*} =
c_{!}\circ (s\times t)^{*}\circ (id\times c)_{!}\circ (id\times (s\times t))^{*}$$
if we knew that the two squares of the above two diagrams were homotopy cartesian.
As this is clearly the associativity condition for the multiplication $\mu$, we see
that it only remains to prove that the two diagrams
$$\xymatrix{X^{(2)} \ar[r]^-{g} \ar[d]_-{f\times (t\circ k)} & X^{(1)} \ar[d]^-{s\times t}\\
X^{(1)}\times X^{(0)} \ar[r]_-{c\times id} & X^{(0)}\times X^{(0)}}$$
$$\xymatrix{
X^{(2)} \ar[r]^-{h} \ar[d]_-{(s\circ f)\times k} & X^{(1)} \ar[d]^-{s\times t}  \\
X^{(0)}\times X^{(1)} \ar[r]_-{id\times c} & X^{(0)}\times X^{(0)}}$$
are homotopy cartesian. This can be easily seen by using
the relations between nerves of categories of equivalences and mapping
spaces in model categories (see for example \cite[Appendix A]{hagII}).
As an example,
let us describe a proof of the fact that the second diagram above is homotopy cartesian (the first one
is even more easy to handle).

We start by a short digression on the notion of fiber product of model categories.
Let
$$\xymatrix{ & M_{1} \ar[d]^-{F_{1}} \\
 M_{2}\ar[r]_-{F_{2}} & M_{3}}$$
be a diagram of model categories and left Quillen functors (we will denote
by $G_{i}$ the right adjoint to $F_{i}$). We define a new model category $M$, called
the \emph{fiber product} of the diagram above and denoted by
$M=M_{1}\times_{M_{3}}M_{2}$, in the following way. The objects in $M$ are
$5$-uples $(x_{1},x_{2},x_{3};u,v)$, where
$x_{i}$ is an object in $M_{i}$, and
$$\xymatrix{ F_{1}(x_{1}) \ar[r]^-{u} &  x_{3} & \ar[l]_-{v} F_{2}(x_{2})}$$
are morphisms in $M_{3}$. Morphisms from $(x_{1},x_{2},x_{3};u,v)$ to
$(x_{1}',x_{2}',x_{3}';u',v')$ in $M$ are defined in the obvious manner,
as families of morphisms $f_{i} : x_{i} \rightarrow x'_{i}$ in $M_{i}$, such that the
 diagram
$$\xymatrix{
F_{1}(x_{1}) \ar[r]^-{u} \ar[d]_-{F_{1}(f_{1})}& x_{3} \ar[d]^-{f_{3}}& F_{2}(x_{2}) \ar[l]_-{v} \ar[d]^-{F_{2}(f_{2})} \\
F_{1}(x_{1}') \ar[r]_-{u'} & x_{3}'& F_{2}(x_{2}') \ar[l]^-{v'}}$$
commutes.

The model structure on $M$ is taken to be the injective levelwise model structure,
for which cofibrations (resp. equivalences) are morphisms $f$ such that
each $f_{i}$ is cofibration (resp. an equivalence) in $M_{i}$ (the fact that this model
structure exists can be checked easily). It is important to note that
an object $(x_{1},x_{2},x_{3};u,v)$ in $M$ is fibrant if and only if
it satisfies the following two conditions
\begin{itemize}
\item The object $x_{i}$ is fibrant in $M_{i}$.

\item The morphisms
$$x_{1} \longrightarrow G_{1}(x_{3}) \qquad x_{2} \longrightarrow G_{2}(x_{3}),$$
adjoint to $u$ and $v$, are fibrations (in $M_{1}$ and $M_{2}$ respectively).

\end{itemize}

Another important fact is that the mapping space of $M$, between two
objects $(x_{1},x_{2},x_{3};u,v)$ and $(x_{1}',x_{2}',x_{3}';u',v')$
is naturally equivalent the following homotopy fiber product
$$\left( Map_{M_{1}}(x_{1},x_{1}')\times Map_{M_{2}}(x_{2},x_{2}') \right)
\times^{h}_{Map_{M_{3}}(\mathbb{L}F_{1}(x_{1}),x_{3}')\times Map_{M_{3}}(\mathbb{L}F_{2}(x_{2}),x_{3}')}
Map_{M_{3}}(x_{3},x_{3}').$$
When the morphisms $u$, $v$, $u'$ and $v'$ are all equivalences
and all objects $x_{i}$, $x_{i}'$ are cofibrant, this homotopy fiber product
simplifies as
$$Map_{M_{1}}(x_{1},x_{1}')\times_{Map_{M_{3}}(x_{3},x_{3}')}^{h}
Map_{M_{2}}(x_{2},x_{2}').$$

We now assume that there exists a
diagram of model categories and left Quillen functors
$$\xymatrix{
M_{0} \ar[r]^-{H_{1}} \ar[d]_-{H_{2}} & M_{1} \ar[d]^-{F_{1}} \\
M_{2} \ar[r]_-{F_{2}} & M_{3},}$$
together with an isomorphism of functors
$$\alpha : F_{1}\circ H_{1} \Rightarrow F_{2}\circ H_{2}.$$
We define a left Quillen functor
$$F : M_{0} \longrightarrow M=M_{1}\times_{M_{3}}M_{2}$$
by sending an object $x\in M$ to $(H_{1}(x),H_{2}(x),F_{2}(H_{2}(x));u,v)$,
where $v$ is the identity of $F_{2}(H_{2}(x))$, and
$$u=\alpha_{x} : F_{1}(H_{1}(x)) \longrightarrow
F_{2}(H_{2}(x))$$
is the isomorphism induced by $\alpha$ evaluated at $x$.

In the situation above, one can restrict to sub-categories of cofibrant objects and
equivalences between them to get a diagram of categories
$$\xymatrix{
wM_{0}^{cof} \ar[r]^-{H_{1}} \ar[d]_-{H_{2}} & wM_{1}^{cof} \ar[d]^-{F_{1}} \\
wM_{2}^{cof} \ar[r]_-{F_{2}} & wM_{3}^{cof},}$$
which is commutative up to the natural isomorphism $\alpha$. Passing to the nerves, one gets
a diagram of simplicial sets
$$\xymatrix{
N(wM_{0}^{cof}) \ar[r]^-{h_{1}} \ar[d]_-{h_{2}} & N(wM_{1}^{cof}) \ar[d]^-{f_{1}} \\
N(wM_{2}^{cof}) \ar[r]_-{f_{2}} & N(wM_{3}^{cof}),}$$
and $\alpha$ induces a natural simplicial homotopy
between $f_{1}\circ h_{1}$ and $f_{2}\circ h_{2}$.

\begin{lem}\label{lstrict}
With the notations above, we assume the following assumption satisfied.
\begin{enumerate}
\item The left derived functor
$$\mathbb{L}F : Ho(M_{0}) \longrightarrow Ho(M)$$
is fully faithful.
\item If an object $(x_{1},x_{2},x_{3};u,v) \in M$ is such that
the $x_{i}$ are cofibrant and the morphisms $u$ and $v$ are
equivalences, then it lies in the essential image of
$\mathbb{L}F$.
\end{enumerate}
Then, the homotopy commutative square
$$\xymatrix{
N(wM_{0}^{cof}) \ar[r]^-{h_{1}} \ar[d]_-{h_{2}} & N(wM_{1}^{cof}) \ar[d]^-{f_{1}} \\
N(wM_{2}^{cof}) \ar[r]_-{f_{2}} & N(wM_{3}^{cof}),}$$
is homotopy cartesian.
\end{lem}

\textit{Proof of the lemma:} This is a very particular case of the strictification theorem
of \cite{sh}, which is reproduced in \cite[Appendix B]{hagII}. We will recall
the proof of the particular case of our lemma
(the general case can be proved in the same
way, but surjectivity on $\pi_{0}$ is more involved).

We consider the induced morphism
$$N(wM_{0}^{cof}) \longrightarrow N(wM_{1}^{cof})\times^{h}_{N(wM_{3}^{cof})}N(wM_{2}^{cof}).$$
In order to prove that this morphism induces a surjective morphism
on $\pi_{0}$ it is well known to be enough to prove the corresponding statement but replacing all the
simplicial sets $N(wM_{i}^{cof})$ by their $1$-truncation. Passing to the fundamental
groupoids, which are the sub-categories of isomorphisms in $Ho(M_{i})$,
we see that it is enough to show that the natural functor of categories
$$Ho(M_{0}) \longrightarrow Ho(M_{1})\times_{Ho(M_{3})}^{h}Ho(M_{2})$$
is essentially surjective. But, an object in the right hand side can be
represented by two cofibrant objects $x_{1} M_{1}^{cof}$,
$x_{2}\in M_{2}^{cof}$, and
an isomorphism $u : F_{1}(x_{1}) \simeq F_{2}(x_{2})$. The isomorphism
$u$ can be itself represented by a string of equivalences in $M_{3}$
$$\xymatrix{ F_{1}(x_{1}) \ar[r] & x_{3} & F_{2}(x_{2}) \ar[l]},$$
where $F_{2}(x_{2}) \rightarrow x_{3}$ is a fibrant model
in $M_{3}$. This data represents an object in $M$, which by  our condition
$(2)$ lifts to an object in $Ho(M_{0})$. This shows that
$$Ho(M_{0}) \longrightarrow Ho(M_{1})\times_{Ho(M_{3})}^{h}Ho(M_{2})$$
is essentially surjective, and thus that
$$N(wM_{0}^{cof}) \longrightarrow N(wM_{1}^{cof})\times^{h}_{N(wM_{3}^{cof})}N(wM_{2}^{cof})$$
is surjective up to homotopy.

It remains to show that the morphism
$$N(wM_{0}^{cof}) \longrightarrow N(wM_{1}^{cof})\times^{h}_{N(wM_{3}^{cof})}N(wM_{2}^{cof})$$
induces an injective morphism in $\pi_{0}$ and isomorphisms on all $\pi_{i}$ for all base
points. Equivalently, it remains to show that for two points $x$ and $y$, the induced
morphism on the path spaces
$$\Omega_{(x,y)}N(wM_{0}^{cof}) \longrightarrow \Omega_{(h_{1}(x),h_{1}(y))}N(wM_{1}^{cof})
\times^{h}_{\Omega_{(f_{1}(h_{1}(x)),f_{1}(h_{1}(y)))}N(wM_{3}^{cof})}
\Omega_{(h_{2}(x),h_{2}(y))}N(wM_{2}^{cof})$$
is a homotopy equivalence. 
Using the relations between path spaces of nerves of equivalences in model
categories and mapping spaces, as recalled e.g. in \cite[Appendix A]{hagII},
we see that the above morphism is equivalent to
$$Map^{eq}_{M_{0}}(x,y) \longrightarrow
Map^{eq}_{M_{1}}(\mathbb{L}F_{1}(x),\mathbb{L}F_{1}(y))\times_{
Map^{eq}_{M_{3}}(\mathbb{L}H_{1}(\mathbb{L}F_{1}(x)),\mathbb{L}H_{1}(\mathbb{L}F_{1}(y)))}^{h}
Map^{eq}_{M_{2}}(\mathbb{L}F_{2}(x),\mathbb{L}F_{2}(y))$$
(here $Map^{eq}$ denotes the sub-simplicial sets of equivalences in the mapping
space $Map$).
It is easy to see that the right hand side is naturally equivalent
to $Map^{eq}_{M}(\mathbb{L}F(x),\mathbb{L}F(y))$, and furthermore that the above morphism
is equivalent to the morphism
$$Map^{eq}_{M_{0}}(x,y) \longrightarrow Map^{eq}_{M}(\mathbb{L}F(x),\mathbb{L}F(y))$$
induced by the left Quillen functor $F$.
The fact that this morphism is an equivalence
simply follows from our condition $(1)$. \hfill $\Box$ \\

We now come back to the proof of the theorem. We need to show that
$$\xymatrix{
X^{(2)} \ar[r]^-{h} \ar[d]_-{(s\circ f)\times k} & X^{(1)} \ar[d]^-{s\times t}  \\
X^{(0)}\times X^{(1)} \ar[r]_-{id\times c} & X^{(0)}\times X^{(0)}}$$
is homotopy cartesian. Simplifying the factor $X^{(0)}$, it is enough to show that
the square
$$\xymatrix{
X^{(2)} \ar[r]^-{h} \ar[d]_-{k} & X^{(1)} \ar[d]^-{t}  \\
X^{(1)} \ar[r]_-{c} & X^{(0)}}$$
is homotopy cartesian. For this we use lemma \ref{lstrict} applied to the
following diagram of model categories and left Quillen functors
$$\xymatrix{
M(T)^{(2)} \ar[r]^-{h} \ar[d]_-{k} & M(T)^{(1)} \ar[d]^-{t} \\
M(T)^{(1)} \ar[r]_-{c} & M(T).}$$
We let $M$ be the fiber product model category
$$M:=M(T)^{(1)} \times_{M(T)}M(T)^{(1)},$$
and let us consider the induced left Quillen functor
$$F : M(T)^{(2)} \longrightarrow M$$
and its right adjoint
$$G : M \longrightarrow M(T)^{(2)}.$$
The objects in the category $M$ are of the form
$(f : a\rightarrow b, g : a' \rightarrow b', a_{0}; u,v)$, where
$f$ and $g$ are morphisms in $M(T)$, $a_{0}$ is an object in $M(T)$, and
$$\xymatrix{b/a \ar[r]^-{u} & a_{0} & b' \ar[l]_-{v}}$$
are morphisms in $M(T)$, where we have denoted symbolically by
$b/a$ the cofiber of the morphism $f$.

The functor $F : M(T)^{(2)} \longrightarrow M$ sends an
object $\xymatrix{x\ar[r] & y \ar[r] & z}$, to
$$(x\rightarrow z, y/x \rightarrow z/x, z/x; id, id)$$
(note that the natural isomorphism $\alpha$ here is trivial).
The right adjoint $G$ to $F$ sends $(a\rightarrow b, a' \rightarrow b', a_{0}; u,v)$
to
$\xymatrix{a\ar[r] & b\times_{a_{0}}a' \ar[r] & b\times_{a_{0}}b'}$.

To finish the proof of theorem \ref{t1} let us check that the left derived functor
$$\mathbb{L}F : Ho(M(T)^{(2)}) \longrightarrow Ho(M)$$
does satisfy the two assumptions of lemma \ref{lstrict}. For this,
we let $A:=\xymatrix{x\ar[r] & y \ar[r] & z}$ be an object in
$M(T)^{(2)}$, and we consider the adjunction morphism
$A \longrightarrow \mathbb{L}G\circ \mathbb{R}F(A)$.
To see that it is an equivalence, it is enough to describe the values
of this morphism levelwise, which are the three natural morphisms in $Ho(M(T))$
$$id : x\rightarrow x \qquad y \longrightarrow 
z\times^{h}_{(z\coprod_{x}^{h}0)}(y\coprod_{x}^{h}0)  \qquad
id : y\rightarrow y.$$
The two identities are of course isomorphisms, and the one in the middle
is an isomorphism because $M(T)$ is a stable model category in the sense
of \cite[\S 7]{ho}.

In the same way, let $B:=(f : a\rightarrow b, g : a' \rightarrow b', a_{0}; u,v)$
be an object in $M$, which is cofibrant (i.e. $f$ and $g$ are cofibrations between
cofibrants objects, and
$a_{0}$ is cofibrant) and such that $u$ and $v$ are equivalences in $M(T)$. We consider
the adjunction morphism in $Ho(M)$
$$\mathbb{L}F\circ \mathbb{R}G(B) \longrightarrow B$$
which we need to prove is an isomorphism. For this we evaluate this morphism
at the various objects $a$, $b$, $a'$, $b'$ and $a_{0}$ to get the following morphisms in $Ho(M(T))$
$$id : a\rightarrow a \qquad id : b\rightarrow b \qquad (b\times_{a_{0}}^{h}a')\coprod_{a}^{h}0 \longrightarrow
a' \qquad (b\times_{a_{0}}^{h}b')\coprod_{a}^{h}0 \rightarrow b' \qquad
(b\times_{a_{0}}^{h}b')\coprod_{a}^{h}0 \rightarrow a_{0}.$$
These morphisms are all isomorphisms in $Ho(M(T))$ as both morphisms $u$ and $v$ are
themselves isomorphisms and moreover because $M(T)$ is a stable model category.
This finishes the proof of the theorem \ref{t1}. \hfill $\Box$ \\

\section{A formula for the derived Hall numbers}

The purpose of this section is to provide a formula
for the numbers $g_{x,y}^{z}$, purely in terms
of the triangulated category $Ho(M(T))$, of $T^{op}$-dg-modules.  \\

Let us assume that $T$ is a locally finite dg-category
over $\mathbb{F}_{q}$, and recall that $Ho(M(T))$
has a natural triangulated structure whose triangles
are the images of the fibration sequences (see \cite[\S 7]{ho}).
We will denote by $[-,-]$ the set of morphisms
in $Ho(M(T))$. The shift functor on $Ho(M(T))$, given by
the suspension, will be denoted by
$$\begin{array}{cccc}
Ho(M(T)) & \longrightarrow & Ho(M(T)) \\
 x & \mapsto & x[1]:=*\coprod^{h}_{x}*.
\end{array}$$
For any two objects $x$ and $y$ in $Ho(M(T))$, we will use the following standard notations
$$Ext^{i}(x,y):=[x,y[i]] \; \forall \; i\in \mathbb{Z}.$$

Let $x$, $y$ and $z$ be three perfect $T^{op}$-dg-modules.
Let us denote by $[x,z]_{y}$ the sub-set of $[x,z]$ consisting of
morphisms $f : x\rightarrow z$ whose cone is isomorphic to $y$ in $Ho(M(T))$.
The group $Aut(x)$, of automorphisms of $x$ in $Ho(M(T))$ acts
on the set $[x,z]_{y}$. For any $f\in [x,z]_{y}$ the stabilizer
of $f$ will be denoted by $Aut(f/z)\subset Aut(x)$. Finally,
we will denote by $[x,z]_{y}/Aut(x)$ the set of orbits of
the action of $Aut(x)$ on $[x,z]_{y}$.

\begin{prop}\label{p3}
With the notation above, one has
$$g_{x,y}^{z}=\sum_{f\in [x,z]_{y}/Aut(x)}|Aut(f/z)|^{-1}.
\prod_{i>0}|Ext^{-i}(x,z)|^{(-1)^{i}}.|Ext^{-i}(x,x)|^{(-1)^{i+1}}$$
$$=\frac{|[x,z]_{y}|.\prod_{i>0}|Ext^{-i}(x,z)|^{(-1)^{i}}}{|Aut(x)|.\prod_{i>0}|Ext^{-i}(x,x)|^{(-1)^{i}}},$$
where as usual $|A|$ denotes the order of the finite set $A$.
\end{prop}

\textit{Proof:} We consider the
morphism
$$c : X^{(1)}(T) \longrightarrow X^{(0)}(T),$$
as well as its homotopy fiber $F^{z}$ at the point
$z$. The simplicial set $F^{z}$ is easily seen to be
naturally equivalent to the nerve of the category
of equivalences in $P(T)/z$, i.e. the nerve of the
category of quasi-isomorphisms between
perfect $T^{op}$-dg-modules over $z$.

We denote by $F^{z}_{x,y}$ the full sub-simplicial set
of $F^{z}$ consisting of all connected components
corresponding to objects $u : x' \rightarrow z$ in $P(T)/z$ such that
$x'$ is quasi-isomorphic to $x$, and the homotopy cofiber of $u$
is quasi-isomorphic to $y$. The simplicial set
$F_{x,y}^{z}$ is locally finite and furthermore
$\pi_{0}(F_{x,y}^{z})$ is finite. By definition one has
$$g_{x,y}^{z}=\sum_{(u : x' \rightarrow z)\in \pi_{0}(F_{x,y}^{z})}
\prod_{i>0}|\pi_{i}(F_{x,y}^{z},u)|^{(-1)^{i}}.$$
It is easily seen that
$\pi_{0}(F_{x,y}^{z})$ is in bijection
with the set $[x,z]_{y}/Aut(x)$. Therefore,
in order to prove the first equality of proposition \ref{p3} it only remains to show that
for a given morphism $u : x \rightarrow z$, considered as
point in $F_{x,y}^{z}$, one has
$$\prod_{i>0}|\pi_{i}(F_{x,y}^{z},u)|^{(-1)^{i}}=
|Aut(f/z)|^{-1}.\prod_{i>0}|Ext^{-i}(x,z)|^{(-1)^{i}}.|Ext^{-i}(x,x)|^{(-1)^{i+1}}.$$
For this, we first notice that there exists homotopy cartesian squares
$$\xymatrix{
Map_{M(T)/z}(x,x) \ar[r] \ar[d] & Map_{M(T)}(x,x) \ar[d] \\
 \bullet \ar[r] & Map_{M(T)}(x,z), }$$
where the bottom horizontal map sends the point $*$ to the
given morphism $u : x \rightarrow z$. Therefore, using the relations between
mapping spaces and nerves of equivalences in model categories
(see \cite[Appendix B]{hagII}), one finds a long exact sequence of finite groups
$$\xymatrix{
Ext^{-i-1}(x,x) \ar[r] & Ext^{-i-1}(x,z) \ar[r] &
\pi_{i+1}(F_{x,y}^{z},u) \ar[r] & Ext^{-i}(x,x) \ar[r] & Ext^{-i}(x,z) \ar[r] &
\dots }$$
$$\xymatrix{ \dots \ar[r] & Ext^{-1}(x,x) \ar[r] & Ext^{-1}(x,z) \ar[r] &
\pi_{1}(F_{x,y}^{z},u) \ar[r] & Aut(f/z) \ar[r] & 0.}$$
This implies the wanted formula
$$\prod_{i>0}|\pi_{i}(F_{x,y}^{z},u)|^{(-1)^{i}}=
|Aut(f/z)|^{-1}.\prod_{i>0}|Ext^{-i}(x,z)|^{(-1)^{i}}.|Ext^{-i}(x,x)|^{(-1)^{i+1}},$$
and thus the first equality of the proposition. Finally, the second equality follows from the general
fact that for a finite group $G$ acting on a finite set $X$, one has
$$\frac{|X|}{|G|}=\sum_{x\in X/G}|G_{x}|^{-1},$$
where $G_{x}$ is the stabilizer of the point $x$.
\hfill $\Box$ \\

For any dg-category $T$ let us denote by $Ho(P(T))$ the full
sub-category of $Ho(M(T))$ consisting of perfect
$T^{op}$-dg-modules. It is a thick sub-category
of $Ho(M(T))$, and thus inherits a natural triangulated structure.

\begin{cor}\label{c2}
Let $T$ and $T'$ be two locally finite dg-categories
over $\mathbb{F}_{q}$. Any triangulated equivalence
$$F : Ho(P(T)) \longrightarrow Ho(P(T'))$$
induces in a functorial way an isomorphism
$$F_{*} : \mathcal{DH}(T) \simeq \mathcal{DH}(T')$$
such that $F_{*}(\chi_{x})=\chi_{F(x)}$.
\end{cor}

\textit{Proof:} Follows from the formula of proposition \ref{p3}.
\hfill $\Box$ \\

Corollary \ref{c2} implies in particular that
the group of self equivalences of the triangulated category
$Ho(P(T))$ acts on the algebra $\mathcal{DH}(T)$.

\begin{rmk}
Corollary \ref{c2} suggests that one can define
the derived Hall algebra of any triangulated category (under
the correct finiteness condition), simply by defining
the derived Hall numbers using the formula of
proposition \ref{p3}.
I do not know if this provides an associative algebra
or not. However, theorem \ref{t1} implies this is the case
for triangulated categories having a dg-model.
\end{rmk}

\section{The heart property}

In this section we examine the behavior of the derived Hall algebra
$\mathcal{DH}(T)$, of a dg-category $T$ endowed with a t-structure, and
the usual Hall algebra of the heart. \\

Let $T$ be a locally finite dg-category. We consider
$Ho(P(T))\subset Ho(M(T))$ the full sub-category of
perfect objects. It is a thick triangulated subcategory
of $Ho(M(T))$. We assume that
$Ho(P(T))$ is endowed with a $t$-structure, compatible with
its triangulated structure, and we let $\mathcal{C}$ be the heart of the
$t$-structure. Recall (e.g. from \cite{de}) that we can construct
a Hall algebra $\mathcal{H}(\mathcal{C})$ of the abelian category $\mathcal{C}$.

\begin{prop}\label{p4}
The Hall algebra $\mathcal{H}(\mathcal{C})$ is naturally isomorphic to
the sub-algebra of $\mathcal{DH}(T)$ spanned by the
$\chi_{x}$ where $x$ lies in the heart $A$.
\end{prop}

\textit{Proof:} We start show that given $x$, $y$ and $z$
three objects in $T$ with $x$ and $y$ in $\mathcal{C}$, then
$$\left( g_{x,y}^{z}\neq 0 \right) \Rightarrow z\in \mathcal{C}.$$
Indeed, if $g_{x,y}^{z}$ is non zero then there exists a triangle
$x\rightarrow z\rightarrow y \rightarrow x[1]$ in $[T]$.
As $x$ and $y$ lie in the heart $\mathcal{C}$ then so does $z$.

It only remains to show that for $x$, $y$ and $z$ in $\mathcal{C}$, our derived Hall numbers
$g_{x,y}^{z}$ coincide with the usual ones. Using 
the explicit  formula Prop. \ref{p3} we find
$$g_{x,y}^{z}=\frac{|[x,z]_{y}|}{|Aut(x)|},$$
which is equal to the number of sub-objects $x'\subset z$ 
such that $x'$ is isomorphic to $x$ and $z/x'$ is isomorphic to 
$y$. \hfill $\Box$ \\

\begin{cor}\label{c1}
Let $A$ be a finite dimensional $k$-algebra, of finite global cohomological
dimension. We let $BA$ be the dg-category having a unique object and
$A$ as its endomorphism ring.
Then, the Hall algebra $\mathcal{H}(A)$ is isomorphic to the
sub-algebra of $\mathcal{DH}(BA^{op})$ generated by the $\chi_{x}$'s, for
$x$ a (non-dg) $A$-module.
\end{cor}

\textit{Proof:} As $A$ is of finite global cohomological dimension,
an $A$-module $M$ is finite dimension as a $k$-vector space
if and only if it is a perfect complex of $A$-modules.
Therefore, the triangulated category $Ho(P(BA^{op}))$ is nothing else
than the bounded derived category
of finite dimensional $A$-modules. We can therefore apply our theorem
\ref{t1} to the standard t-structure on $Ho(P(BA^{op}))$, whose heart is
precisely the abelian category of finite dimensional $A$-modules. \hfill $\Box$ \\

An important example of application of theorem \ref{t1} is
when $T$ is the dg-category of perfect complexes on
some smooth projective variety $X$ over $k=\mathbb{F}_{q}$.
This provides an associative algebra $\mathcal{DH}(X):=\mathcal{DH}(T)$,
containing the Hall algebra of the category of coherent sheaves on $X$.

\section{Derived Hall algebra of hereditary abelian categories}

In this last section we will describe the derived Hall algebra
of an hereditary (i.e. of cohomological dimension $\leq 1$) abelian category $A$, 
in terms of its underived Hall algebra. \\

We let $T$ be a locally finite dg-category over $\mathbb{F}_{q}$, and
we assume that we are given a t-structure on the triangulated
category $Ho(P(T))$, as in the previous section. We let 
$\mathcal{C}$ be the heart of this t-structure, and we will 
assume that $\mathcal{C}$ is an hereditary abelian category, or
in other words that for any $x$ and $y$ in $\mathcal{C}$ we have
$Ext^{i}(x,y)=0$ for all $i>1$ (the ext groups here are
computed in $\mathcal{C}$). We will also assume that the natural
triangulated functor $D^{b}(\mathcal{C}) \longrightarrow Ho(P(T))$
(see \cite{bbd}) is an equivalence, or equivalently that the t-structure
is bounded and that for any 
$x$ and $y$ in the heart we have $[x,y[i]]=0$ for any $i>1$ 
(here $[-,-]$ denotes the set of morphisms in $Ho(P(T))$). 
According to \cite{de} we can construct 
a Hall algebra $\mathcal{H}(\mathcal{C})$, and our purpose is to 
describe the full derived Hall algebra $\mathcal{DH}(T)$ by generators
and relations in terms of $\mathcal{H}(\mathcal{C})$. 

For this, we need to introduce the following notations. For two 
objects $x$ and $y$ in $\mathcal{C}$, we set as usual

$$\left\langle x,y \right\rangle :=Dim_{\mathbb{F}_{q}}Hom(x,y)-
Dim_{\mathbb{F}_{q}}Ext^{1}(x,y).$$

For four objects $x$, $y$, $c$ and $k$ in $\mathcal{C}$, we 
define a rational number $\gamma_{x,y}^{k,c}$ as follows. 
We first consider $V(k,y,x,c)$, the subset of 
$Hom(k,y)\times Hom(y,x) \times Hom(x,c)$ consisting of all
exact sequences
$$\xymatrix{ 0 \ar[r] & k \ar[r] & y \ar[r] & x \ar[r] & c \ar[r] & 0.}$$
The set $V(k,y,x,c)$ is finite and we set 
$$\gamma_{x,y}^{k,c}:=\frac{|V(k,y,x,c)|}{|Aut(x)|.|Aut(y)|}.$$

Recall finally that for three objects $x$, $y$ and $z$ in $\mathcal{C}$
we define $g_{x,y}^{z}$ as being the number of sub-objects
in $x' \subset z$ such that $x'$ is isomorphic to $x$ and
$z/x'$ is isomorphic to $y$. By our proposition \ref{p4}, 
this number $g_{x,y}^{z}$ is also the derived
Hall number of the corresponding $T^{op}$-perfect dg-modules. 

\begin{prop}\label{p5}
The $\mathbb{Q}$-algebra $\mathcal{DH}(T)$ is isomorphic
to the associative and unital algebra generated by the set
$$\{Z_{x}^{[n]}\}_{x\in \pi_{0}(\mathcal{C}),n\in \mathbb{Z}},$$
where $\pi_{0}(\mathcal{C})$ denotes the set of isomorphism 
classes of objects in $\mathcal{C}$, 
and with the following relations
\begin{enumerate}
\item $$Z_{x}^{[n]}.Z_{y}^{[n]}=\sum_{z  \in \pi_{0}(\mathcal{C})}g_{x,y}^{z}.Z_{z}^{[n]}$$

\item $$Z_{x}^{[n]}.Z_{y}^{[n+1]}=\sum_{c,k \in \pi_{0}(\mathcal{C})}
\gamma_{x,y}^{k,c}.q^{-\left\langle c,k \right\rangle}.Z_{k}^{[n+1]}.Z_{c}^{[n]}$$

\item For $n-m<-1$

$$Z_{x}^{[n]}.Z_{y}^{[m]}=q^{(-1)^{n-m}.\left\langle x,y \right\rangle}.Z_{y}^{[m]}.Z_{x}^{[n]}.$$

\end{enumerate}
\end{prop}

\textit{Proof:} We let $\pi_{0}(P(T))$ be the set of isomorphism
classes of objects in $Ho(P(T))$, and for $x\in \pi_{0}(P(T))$
we let $\chi_{x} \in \mathcal{DH}(T)$ the characteristic function
of $x$. In the same way, for $x\in \pi_{0}(\mathcal{C})$ we denote
by $\chi_{x} \in \mathcal{H}(\mathcal{C})$ the corresponding element. 
By our assumption on $T$, any object in 
$Ho(P(T))$ can be written as a direct sum
$\oplus_{n}x_{n}[n]$, where 
$x_{n}$ lies in the heart $\mathcal{C}$.

We construct a morphism of vector spaces
$$\Theta : B \longrightarrow \mathcal{DH}(T),$$
where $B$ is the algebra generated by the $Z_{x}^{[n]}$
as described in the proposition, by setting
$$\Theta(Z_{x}^{[n]}):=\chi_{x[n]}.$$
We claim that $\Theta$ is an isomorphism of algebras. 

To prove that $\Theta$ is a morphism of algebras we need
to check that the corresponding relations in 
$\mathcal{DH}(T)$ are satisfied. In other words we have to prove the
following three relations are satisfied

\begin{enumerate}
\item $$\chi_{x[n]}.\chi_{y[n]}=\sum_{z  \in \pi_{0}(\mathcal{C})}g_{x,y}^{z}.\chi_{z[n]}$$

\item $$\chi_{x[n]}.\chi_{y[n+1]}=\sum_{c,k \in \pi_{0}(\mathcal{C})}
\gamma_{x,y}^{k,c}.q^{-\left\langle c,k \right\rangle}.\chi_{k[n+1]}.\chi_{c[n]}$$

\item For $n-m<-1$

$$\chi_{x[n]}.\chi_{y[m]}=q^{(-1)^{n-m}.\left\langle x,y \right\rangle}.\chi_{y[m]}.\chi_{x[n]}.$$

\end{enumerate}

The relation 
$(1)$ follows from the heart property Prop. \ref{p4}. For
relation $(3)$, the only extension of $y[m]$ by $x[n]$
is the trivial one, as $[y[m],x[n+1]]=0$ because 
$n-m<-1$. In the same way, we have $[x[n],y[m+1]]=0$, because
$\mathcal{C}$ is of cohomological dimension $\leq 1$. 
Therefore, in order to check the relation $(3)$ we only need to check that 
$$g_{x[n],y[m]}^{x[n]\oplus y[m]}=q^{(-1)^{n-m}.\left\langle x,y \right\rangle} \qquad 
g_{y[m],x[n]}^{x[n]\oplus y[m]}=1,$$
which is a direct consequence of the formula given in Prop. \ref{p3}.

It remains to prove relation $(2)$, which can be easily be rewritten as
$$\chi_{x[n]}.\chi_{y[n+1]}=\sum_{c,k \in \pi_{0}(\mathcal{C})}
\gamma_{x,y}^{k,c}.q^{-\left\langle c,k \right\rangle}.\chi_{k[n+1]\oplus c[n]}.$$
By shifting if necessary it is in fact enough to prove that
$$\chi_{x}.\chi_{y[1]}=\sum_{c,k \in \pi_{0}(\mathcal{C})}
\gamma_{x,y}^{k,c}.q^{-\left\langle c,k \right\rangle}.\chi_{k[1]\oplus c}.$$
An extension of $y[1]$ by $x$ is classified by a morphism
$y[1] \rightarrow x[1]$, and therefore by a morphism
$u : y \rightarrow x$. For such a morphism $u$, we let 
$c$ and $k$ be defined by the long exact sequence in $\mathcal{C}$
$$\xymatrix{
0 \ar[r] & k \ar[r] & y \ar[r]^-{u} & x \ar[r] & c \ar[r] & 0.}$$
Clearly, for a morphism $u : y \rightarrow x$, and the corresponding
triangle $\xymatrix{x \ar[r] & z \ar[r] & y[1] \ar[r]^-{u[1]} & x[1]}$, 
the object $z$ is isomorphic in $Ho(P(T))$ to 
$k[1]\oplus c$. Therefore, we can write
$$\chi_{x}.\chi_{y[1]}=\sum_{c,k \in \pi_{0}(\mathcal{C})}
g_{x,y[1]}^{k[1]\oplus c}.\chi_{k[1]\oplus c},$$
and it only remains to show that 
$$g_{x,y[1]}^{k[1]\oplus c}=\gamma_{x,y}^{k,c}.q^{-\left\langle c,k \right\rangle}.$$
Applying Prop. \ref{p3} we have
$$g_{x,y[1]}^{k[1]\oplus c}=
\frac{|[x,k[1]\oplus c]_{y}|}{|Aut(x)|.|[x,k]|}.$$
The set $[x,k[1]\oplus c]_{y}$ can be described 
as the set of ordered pairs $(\alpha,\beta)$, where 
$\alpha : x \rightarrow c$ is an epimorphism in $\mathcal{C}$, 
$\beta \in [x,k[1]]$ is an extension 
$\xymatrix{0 \ar[r] & k \ar[r] & y' \ar[r] & x \ar[r] & 0}$, 
such that the object $y$ is isomorphic to $y'\times_{x} ker(\alpha)$. 
For each epimorphism $\alpha : x \rightarrow c$, 
we form the short exact sequence
$$\xymatrix{
0 \ar[r] & E_{\alpha} \ar[r] & [x,k[1]] \ar[r] & [ker(\alpha),k[1]] \ar[r] & 0.}$$
We see that the set $[x,k[1]\oplus c]_{y}$ is in bijection with the 
set of triples $(\alpha,\beta',a)$, where 
$\alpha : x \rightarrow c$ is an epimorphism, 
$\beta'$ is an extension class 
$\xymatrix{0 \ar[r] & k \ar[r] & y' \ar[r] & ker(\alpha) \ar[r] & 0}$
such that $y'$ is isomorphic to $y$, and $a$ is an element 
in $E_{\alpha}$. 

Let us now fix an epimorphism $\alpha : x \rightarrow c$, and
let $V(k,y,\ker(\alpha))$ be the set of exact sequences
$$\xymatrix{
0 \ar[r] & k \ar[r] & y \ar[r] & ker(\alpha) \ar[r] & 0.}$$
The number of possible $\beta'$ is then equal to 
$$\frac{|V(k,y,\ker(\alpha))|}{|Aut(y)|}.|[ker(\alpha),k]|,$$
as $[ker(\alpha),k]$ is (non-canonically) the stabilizer of
any point in $V(k,y,\ker(\alpha))$ for the action of $Aut(y)$. 
Furthermore, considering the long exact sequence
$$\xymatrix{
0 \ar[r] & [c,k] \ar[r] & [x,k] \ar[r] & [ker(\alpha),k] \ar[r] & 
[c,k[1]] \ar[r] & [x,k[1]] \ar[r] & [ker(\alpha),k[1]] \ar[r] & 0}$$
we see that the number of points in the set $E_{\alpha}$
is equal to
$$|E_{\alpha}|=|[c,k[1]]|.|[ker(\alpha),k]|^{-1}.
|[x,k]|.|[c,k]|^{-1}.$$
All together, this implies that 
$$g_{x,y[1]}^{k[1]\oplus c}=\frac{|[x,k[1]\oplus c]_{y}|}{|Aut(x)|.|[x,k]|}=
\sum_{\alpha : x \twoheadrightarrow c}
\frac{|V(k,y,\ker(\alpha))|}{|Aut(y)|.|Aut(x)|}.|[ker(\alpha),k]|.|[x,k]|^{-1}.|E_{\alpha}|$$
$$=
\left( \sum_{\alpha : x \twoheadrightarrow c}
\frac{|V(k,y,\ker(\alpha))|}{|Aut(y)|.|Aut(x)|}\right).|[c,k[1]]|.|[c,k]|^{-1}
=\gamma_{x,y}^{k,c}.q^{-\left\langle c,k \right\rangle}.$$
This finishes the proof of the fact that $\Theta$ is a morphism 
of algebras. 

To prove that $\Theta$ is
bijective we consider the following diagram
$$\xymatrix{
\bigotimes_{\mathbb{Z}}\mathcal{H}(\mathcal{C}) \ar[r]^-{u} \ar[d]_-{v} & B \ar[dl]^-{\Theta} \\
\mathcal{DH}(T), & }$$
where the infinite tensor product is the restricted one (i.e. 
generated by words of finite lengh). The morphisms
$u$ and $v$ are defined by
$$u(\otimes_{n} \chi_{x_{n}}):=\prod^{\rightarrow}_{n} Z_{x_{n}}^{[n]}$$
$$v(\otimes_{n} \chi_{x_{n}}):=\chi_{\oplus_{n}x_{n}[n]},$$
where 
$\prod^{\rightarrow}_{n}Z_{x_{n}}^{[n]}$ means the ordered 
product of the elements $Z_{x_{n}}^{[n]}$
$$\prod^{\rightarrow}_{n}Z_{x_{n}}^{[n]}=\dots
Z_{x_{i}}^{[i]}.Z_{x_{i-1}}^{[i-1]}.\dots.$$

For any two objects $x$ and $y$ in $Ho(P(T))$
with $x\in Ho(P(T))^{\leq 0}$ and $y\in Ho(P(T))^{>0}$ we have
$$g_{x,y}^{x\oplus y}=1 \qquad g_{x,y}^{z}=0 \quad if \; \chi_{z}\neq \chi_{x\oplus y}.$$
In other words $\chi_{x}.\chi_{y}=\chi_{x\oplus y}$. This implies on the one hand that 
the diagram above commutes, and on the other hand
that the morphism $v$ is an isomorphism. It remains to show that 
the morphism $u$ is also an isomorphism, but surjectivity is clear
and injectivity follows from injectivity of $v$.  \hfill $\Box$ \\

It seems to us interesting to note the analogies between the
description of $\mathcal{DH}(T)$ given by Prop. \ref{p5}
and the Lattice algebra $L(T)$ of Kapranov introduced in \cite{ka}. 
The formal structures of the two algebras $\mathcal{DH}(T)$ and $L(T)$ seem
very close, thought $L(T)$ contains an extra factor corresponding the
\emph{Cartan subalgebra} (i.e. a copy of the group 
algebra over the Grothendieck group of $T$). Therefore, it seems likely that 
there exists an embedding $j : \mathcal{DH}(T)\hookrightarrow L(T)$. 
However, I do not see any obvious way to construct 
such a embedding, and I prefer to leave the question of comparing 
$\mathcal{DH}(T)$ and $L(T)$ for further research.

\end{document}